\documentclass[12pt]{article}
\usepackage{amsfonts}
\usepackage{graphicx}   
\usepackage{amsmath}

\begin{document}

\begin{center}
{\Large \textbf{Elementary trigonometry based on a first order differential equation}} 

Horia I. Petrache  

Department of Physics, Indiana University Purdue University Indianapolis, Indianapolis, IN 46202

\today

\end{center}

\textbf{Abstract}

It is shown that with appropriate boundary conditions, a real function satisfying the differential equation $f'(x) = f(x+a)$ has all known properties of the sine function. A number of elementary derivations are presented including proofs for periodicity and trigonometric identities. 

\section{Introduction}

In elementary analysis, the sine and cosine functions are identified as the two linearly independent solutions of the second order differential equation $f'' = -f$. Therefore, any function $f$ that satisfies this harmonic equation is periodic, bounded, and satisfies a number of trigonometric identities. However, many of these trigonometric properties are derived not from $f'' = -f$ but from alternative definitions of sine and cosine functions. These alternatives include primarily geometrical definitions based on right triangles and unit circles [1,2,3], power series and complex exponentials [4,5], and integral forms of arc lengths [6,7]. The question then is what properties follow directly from $f'' = -f$. For example, it is quite easy to show by direct differentiation that $f^2+f'^2 = const.$, which then implies that $f$ and $f'$ are both bounded. However, other properties are more difficult to prove in abstract form, e.g. the fact that $f$ is periodic. An even more interesting question is whether there is any other property of the sine and cosine functions from which all known properties can be derived. Here it is shown that the differential equation $f'(x) = f(x+a)$ specifies the sine function given the appropriate boundary conditions. Derivations are carried at an abstract level that does not involve series expansions or complex exponentials.

\section{Derivations}

Consider the differential equation
\begin{equation}
f'(x) = f(x+a),
\label{eq:fxa}
\end{equation}
in which the argument $x$ and the constant $a$ are real. In the following treatment we assume that $f$ is well defined on the entire real axis, is of class $C^2$ or higher, and it is not identically zero. 
 
Let us first note that the differential operator changes the parity (symmetry) of an odd or even function. If $f$ is odd then $f'$ is even and vice versa. For the above equation to hold, the translation operator that takes $f(x)$ into $f(x+a)$ must perform the same symmetry transformation as the differential operator. From this perspective, we will investigate solutions of Eq.~\ref{eq:fxa} with well defined parity. 

\textbf{Proposition 1.} \textit{$f'(x) = f(x+a)$ implies that $f$ is periodic.}

\textbf{Proof.} Use notation $g(x) = f'(x)$ and assume that $f$ is odd. (The proof for an even function is similar.) If $f$ is odd then $g=f'$ is even. We have $f(x) \overset{Eq.~\ref{eq:fxa}}{=} f'(x-a) \overset{def}{=} g(x-a) \overset{sym}{=} g(-x+a) \overset{Eq.~\ref{eq:fxa}}{=} f(-x+2a) \overset{sym}{=} -f(x-2a)$, where the text above the equal sign indicates the property used at each step, and where $sym$ stands for symmetry and $def$ for the definition (notation) $g = f'$. We have obtained that $f(x) = -f(x-2a)$ for any $x$, which implies that 
$f(x-2a) = -f(x-4a)$, so that $f(x) = f(x-4a)$ for any $x$. It means that $f$ is periodic of period $4a$.

\textbf{Proposition 2.} \textit{$f'(x) = f(x+a)$ implies that $g'(x) = g(x+a)$, where $g = f'$.}

\textbf{Proof.} $g'(x) = f''(x) \overset{Eq.~\ref{eq:fxa}}{=} f'(x+a) = g(x+a)$.

\textbf{Proposition 3.} \textit{$f'(x) = f(x+a)$ implies that $g'(x) = -f(x)$, where $g = f'$.}

\textbf{Proof.} Using Proposition 2, we have 
$g'(x) = g(x+a) \overset{sym}{=} g(-x-a) \overset{def}{=} f'(-x-a) \overset{Eq.~\ref{eq:fxa}}{=} f(-x-a+a) = f(-x) \overset{sym}{=} -f(x).$ We have in fact proven that $f$ satisfies the second order differential equation
\begin{equation}
f'' = -f.
\label{eq:h}
\end{equation}

\textbf{Proposition 4.} \textit{$f^2 +g^2 = constant$.}

Need to show that the derivative of $f^2 +g^2$ is zero. We have 
$ff' + gg' = ff' + f'(-f) = 0$. Because the functions $f$ and $g$ are specified up to an arbitrary scale, we can choose the constant to be 1 without loss of generality. If $f$ is odd, we have $f(0) = 0$ and therefore $g(0) = \pm 1$ and we can choose $g(0) = 1$.

\textbf{Proposition 5.} \textit{The functions $f$ and $g=f'$ are linearly independent.}

\textbf{Proof.} Assume that $f = cf'$ with $c$ a real constant. We have $f' = cf'' = -cf = -c^2f'$, implying that $c^2 = -1$ which contradicts the assumption that $c$ was real. 
(Alternatively, it can be verified that the Wronskian of $f$ and $f'$ is non-zero.)

\textbf{Proposition 6.}  \textit{If $f$ is a solution of $f'' = -f$, then $g = f'$ is also a solution and $f$ and $g$ are linearly independent.} 

\textbf{Proof.} $g'' = f''' = -f' = -g$. Linear independence was proven above.

\textbf{Proposition 7.}  \textit{$f'' = -f$ implies that $f(x+y) = f'(y)f(x)+f(y)f'(x)$.}

\textbf{Proof.} First note that $f(x+y)$ is a solution of $f'' = -f$ for any $y$, where the derivative is with respect to the variable $x$. That means that $f(x+y)$ is a linear combination of $f$ and $f'$, which were shown earlier to be independent solutions of Eq.~\ref{eq:h}. We have
\begin{equation}
f(x+y) = A(y)f(x) + B(y)f'(x),
\end{equation}
where the coefficients $A$ and $B$ are functions of $y$. With $f(0) = 0$ and $f'(0) = 1$, we have 
\begin{eqnarray}
f(0+y) &=& A(y)f(0) + B(y)f'(0) \\
f(x-x) & = & A(-x)f(x) + B(-x)f'(x),
\end{eqnarray}
which gives $B(y) = f(y)$ and $A(-x)f(x) + B(-x)f'(x) = 0$. Therefore, 
$A(-x)f(x) -f(x)f'(x) = 0$ for any $x$, so $A(-x)-f'(x) = 0$, giving $A(x) = f'(-x) = f'(x)$.

For any $x$ and $y$ we have obtained that
\begin{equation}
f(x+y) = A(y)f(x) + B(y)f'(x) = f'(y)f(x)+f(y)f'(x),
\end{equation}

\textbf{Proposition 8.}  $f'' = -f$ \textit{implies that} $f_{\mbox{max }}f'(x) = f(x+a)$, \textit{where} $f_{\mbox{max}} = f(a)$.

\textbf{Proof.} If $f(a) = f_{\mbox{max }}$ then $f'(a) = 0$ (Proposition 4) and using Proposition 7, 
$f(x+a) = f'(a)f(x)+f(a)f'(x) = f_{\mbox{max }} f'(x)$. Choosing $f_{\mbox{max}} =1$, we obtain $f'(x) = f(x+a)$.

\section{Conclusions}

Based on the above development, we can conclude that the differential equations $f'(x) = f(x+a)$ and $f'' = -f$ are equivalent in that each can be regarded as a fundamental property that defines the sine function. To be more precise, each equation defines a general linear combination of sine and cosine from which the sine function is selected with the appropriate choice of symmetry and overall amplitude. The differential equations satisfied by the sine and cosine functions have theoretical significance not only in mathematics but also in physics and other applied fields. However, as opposed to the well known harmonic equation that describes the motion of oscillators, Eq.~\ref{eq:fxa} appears under-appreciated. As presented here, Eq.~\ref{eq:fxa} represents the particular case in which the action of two distinct operators (differential and translation operator) coincides on a subset of state functions. 

Acknowledgments

The author thanks Pat DeMoss, Heather Stout, and Kashyap Vasavada for valuable discussions and suggestions for the manuscript. 

\textbf{References}

[1] A. R. Crathorne and E. B. Lytle, Trigonometry,  Henry Holt and Company, New York, 1938. 

[2] E. P. Vance, Trigonometry, Addison-Wesley Publishing, Cambridge, 1954. 

[3] H. I. Petrache, Correspondence between geometrical and differential definitions of the sine and cosine functions and connection with kinematics, arXiv:1110.6106v1 [physics.ed-ph]. 

[4] W. Rudin, Principles of Mathematical Analysis, McGraw-Hill Inc., New York, 1964.

[5] J. F Randolph,  Basic real and abstract analysis,  Academic Press, New York, 1968.

[6] K. R. Stromberg, An Introduction to Classical Real Analysis, Waldsworth Inc., Belmont, CA, 1981.

[7] J. M. H. Olmsted,  Real variables, an introduction to the theory of functions, Appleton-Century-Crofts, New York, 1959.
 
\end{document}